\documentclass[12pt,a4paper]{article}
\usepackage[a4paper, left=2.1cm]{geometry}
\usepackage{amsfonts}
\usepackage{amsmath}
\usepackage{amssymb}
\usepackage{amsthm}
\usepackage[utf8]{inputenc}
\date{}
\title {Continuants with equal values, a combinatorial approach}

\author
{G Ramharter, L.Q. Zamboni}

\begin{document}
\pagenumbering{arabic}
\maketitle

\begin{abstract}
{\small A regular  continuant is the denominator $K$ of a terminating regular continued fraction, interpreted as a function of the partial quotients. We regard $K$ as a function defined on the set of all finite words on the alphabet $1<2<3<\dots$ with values in the positive integers.
Given a word $w=w_1\cdots w_n$ with $w_i\in\mathbb{N}$ we define its multiplicity $\mu(w)$  as the number of times the value $K(w)$ is assumed in the Abelian class $\mathcal{X}(w)$ of all permutations of the word $w.$ 
We prove that there is an infinity of different lacunary alphabets of the form 
$\{b_1<\dots <b_t<l+1<l+2<\dots <s\}$ with $b_j, t, l, s\in\mathbb{N}$ and $s$ sufficiently large such that $\mu$ takes arbitrarily large values for words on these alphabets. 
The method of proof relies in part on a combinatorial characterisation of the word $w_{max}$ in the class $\mathcal{X}(w)$ where $K$ assumes its maximum. 

\medskip
\parindent0.0cm
MSC: primary 11J70; secondary 68R15, 68W05, 05A20.

}

\end{abstract}

\smallskip
\parindent0.0cm
\vskip-0.2cm                                                                                   
\parindent0cm
\bigskip
{\bf Introduction.} Given a sequence $w=(w_1,\ldots , w_n),$ of positive $w_i,$  let $K(w)$ be the continuant of $w,$ 
i.e., the denominator of the finite regular continued fraction 
$\frac{1}{w_1+}\frac{1}{w_2+}\dots\frac{1}{w_{n-1}+} \frac{1}{w_n}$.
We shall regard $w$ as a {\it word} of length $n$ over the alphabet $\{1<2<3<\dots\}$ and write $w = w_1\cdots w_n.$
Since $K(w)=K(\bar{w}),$ where $\bar{w}=w_n\cdots w_1$ denotes the reversal of $w,$ we shall henceforth identify each word $w$ with its reverse 
$\bar{w}$. 
Let $\mathcal{X}(w)$ denote the Abelian class of $w$ consisting of all permutations of $w.$
The following problem has attracted much attention and led to a number of applications (see e.g.\ \cite{Baxa, Ramh_Some, Ramh_FineWilf, ShaSor, Zam}): 
Let $A = \{a_1<\dots < a_s\}$ be a finite ordered alphabet with $a_j\in\mathbb{N}$.
Given a word $w=w_1w_2\cdots w_n$ with $w_i\in A$, find the arrangements $w_{\max}, w_{\min}\in\mathcal{X}(w)$ maximizing resp.\ minimizing the function $K(\cdot)$ on $\mathcal{X}(w)$. 
The first author \cite{Ramh_Extr} gave an explicit description of both extremal arrangements $w_{\max}$ and $w_{\min}$ and showed that in each case the arrangement is unique (up to reversal) and independent of the actual values of the positive integers $a_i.$ 
He also investigated the analogous problem for the semi-regular continuant $K'$ defined as the denominator of the semi-regular continued fraction $\frac{J'}{K'}=\frac{1}{w_1-}\frac{1}{w_2-}\dots\frac{1}{w_{n-1}-}\frac{1}{w_n}$ with entries $w_i\in \{2, 3, ...\}.$
He gave a fully combinatorial description of the minimizing arrangement $w'_{\min}$ for $K'(\cdot)$ on $\mathcal{X}(w)$ and  showed that the arrangement is unique (up to reversal) and independent of the actual values of the positive integers $a_i.$
However, the determination of the maximizing arrangement $w'_{\max}$ for the semi-regular continuant turned out to be more difficult.
He showed that in the special case of a $2$-digit alphabet $\{(2\leq )~a_1<a_2\},$ the maximizing arrangement $w'_{\max}$ is a Sturmian word and is independent of the values of the $a_i.$ 
Recently the second author together with M. Edson and A. De Luca \cite{Zam} developed an algorithm for constructing $w'_{\max}$ over any ternary alphabet $\{(2\leq )~a_1 <a_2 < a_3\}$, and  showed that the maximizing arrangement is independent of the choice of the digits.
In contrast,  they exhibited examples of words $w=w_1\cdots w_n$ over a $4$-digit alphabet $A=\{(2\leq )~a_1<a_2<a_3<a_4\}$ for which the maximizing arrangement for $K'(\cdot)$ is not unique and depends on the actual values of the positive integers $a_1$ through $a_4.$ 
In the course of these investigations the following problem came up: 
given an alphabet $A$ of positive integers, we say that a word $w$ on $A$ has multiplicity $\mu=\mu (w)$ if the value $K(w)$ occurs precisely $\mu$ times in the multi-set $\{ K(x): x\in\mathcal{X}(w)\}$. The multiplicity $\mu' (w)$ is defined analogously for the semi-regular continuant $K'(w).$ 
Thus each Abelian class $\mathcal{X}(w)$ is split into subclasses of equally valued words. 
Question: is it true that $\mu$ can take arbitrarily large values for infinitely many alphabets and is there a combinatorial proof of this? 
Our aim here is to give a positive answer to this question in the case of regular continued fractions.

\vskip0.45cm
{\bf Theorem}.                                                             
Fix positive integers 
\hskip0.3mm 
$1\leq t\hskip-0.5mm \leq\hskip0mm l\hskip-0.5mm < \hskip-0.5mm s,\hskip1mm b_1\hskip-1mm <\hskip-0.7mm ...\hskip-0.7mm <\hskip-0.6mm b_t\hskip-0.5mm \leq\hskip-0.5mm l$ and let $A$ be an ordered alphabet of the form 
$\{b_1<...<b_t<l+1\hskip-0.5mm <...\hskip-0.2mm <\hskip-0.5mm s\}$.
Then for all $s$ sufficiently large, there exists an infinite sequence of words $w_k$ over $A$ with multiplicities 
$\mu (w_k)\to\infty$ as $k\to\infty$.

\bigskip
It should be noted that for fixed $s$ one obtains the largest possible alphabet $A'\hskip-0.8mm=\hskip-0.8mm\{1\hskip-0.7mm<\hskip-0.7mm 2<\hskip-1mm\dots< s\}$ by choosing $b_1=t=l=1\ (< s)$. 
Our proof makes use of the combinatorial structure of $w_{\max}$ found by the first author in \cite{Ramh_Extr}. 

\bigskip
{\bf Preliminaries.} We introduce some notation. Let $w=w_1\cdots w_n$ be a word of length $n\geq 2$ 
with $w_j\in\mathbb{N}$  ($j = 1,\dots, n)$. 
The regular continuant of $w$ has a matrix representation 
 \[K(w_1) = w_1\ \ \mbox{and}\ \ 
K(w)=\det\left( \begin{matrix}
w_1 & -1 & 0 &\cdots &0           \\
1 & w_2 & -1 &\ddots&\vdots       \\
0 & 1 & \ddots &\ddots &0         \\
\vdots &\ddots&\ddots &w_{n-1}& -1\\
0 &\cdots&0&1&w_n
\end{matrix} \right),\ n\geq 2\] 

It can also be defined recursively by 
$K(\{\hskip0.05cm\}) = 1$ ($\{\hskip0.05cm\}$ = empty word), $K(w_1) = w_1$ and 
$K(w_1 \cdots w_j) = w_{j} K(w_1 \cdots w_{j-1})+K(w_1 \cdots w_{j-2})$ for $j\geq 2$.
For each $1\leq k \leq m\leq n$ we set 
$w_{k,m}:=w_k\cdots w_m$ and $W:=W_{1,n},~W_{k,m} := K(w_{k,m}).$ 
The following fundamental formula  goes back to the late 19$^{{\rm th}}$ century and can be found in 
Perron \cite{Perron}, p.11, (4): 
$(W\hskip-0.8mm =\hskip0.5mm)~W_{1,n}=W_{1,j}W_{j+1,n}+W_{1,j-1}W_{j+2,n}$ 
$(j\in\{1,\dots ,n-1\}).$ 
From this we infer the simple but useful inequality
\vskip-0.20cm
\begin{equation}\label{ineq1}    
W_{1,n}<2W_{1,j}W_{j+1,n}.
\end{equation}
\vskip0.1cm
Let $A = \{a_1<\dots <a_{s}\}\subset \mathbb{N}.$ 
We consider a word $w = w_1\cdots w_n := a_1^{p_1}\hskip-0.4mm \cdots a_{s}^{p_{s}}$ 
of length $n$ with Parikh vector ${\bf p} = (p_1,...,p_s)$ with  $p_1+\cdots +p_s = n$ 
where 
\vskip-0.15cm
\[a^r=\underbrace{aa\cdots a}_{\mbox{$r$-times}}\]
\vskip0.15cm
denotes a sequence of $r$ equal elements $a$. Let $\mathcal{X}=\mathcal{X}(A,{\bf p})$ 
denote the set of all permutations of 
$w$ where we identify each word $v$ with its reverse $\bar v.$ Let $N(A,{\bf p})$  denote the cardinality of $\mathcal X.$ 
Then, $N(A,{\bf p}) \geq \frac{n!}{2p_1!\dots p_s!}.$ 
We put $W_{\max}=W_{\max}(A,{\bf p}):=\max\{K(v): v\in\mathcal{X}\}$. 
It was shown in \cite{Ramh_Extr} (see (3), p.~190)  that $W_{\max}$ is uniquely attained (up to reversal) by the arrangement 
\vskip-0.1cm
\begin{equation}\label{eqn2} %
a_s L_{s-1}a_{s-2}L_{s-3}\cdots a_1^{p_1}\cdots  a_{s-3}L_{s-2}a_{s-1}L_s\end{equation}
\vskip0.2cm
where $L_i=a_i^{p_i-1}.$
Let $P=P(A,{\bf p})=\#\{K(v) : v\in \mathcal{X}\}.$
\vskip0.4cm 
\begin{proof}
[{\bf Proof of the Theorem}]                                             
Our first goal is to describe how to specify the last digit $s\ (\geq 2)$ in an alphabet 
$A : \{b_1<\dots <b_t <l+1< \cdots <s\}$. 
We consider  'equipartitioned' words
\vskip0.cm
\begin{displaymath}
w=w_1\cdots w_n := b_1^m\cdots b_t^m (l+1)^m\cdots s^m.
\end{displaymath}
\vskip0.3cm
corresponding to the Parikh vector ${\bf p}= (m,m,\ldots,m)$
in which each digit of $A$ occurs precisely $m$-times in $w.$ 
We will give a lower  bound for $s$ (see (\ref{ineq7}) below).\  
To this end, we introduce the quantities 
$Q_{r,m-1}:=K(r^{m-1})~(r\in {1, 2, ...})\hskip0.5mm$. 
They are the elements of the $r$-th generalised Fibonacci sequence which is determined by the recursion 
$Q_{r,0}:= 1$, $Q_{r,1}:=K(r)=r$, $Q_{r,j+1}:=rQ_{r,j}+Q_{r,j-1}\ (j=1, 2, \dots)$. \\

{\bf Claim:} $Q_{r,j-1} < (r+1)^{j}$ for each fixed $r\geq 1$ and all $j\geq 1.$  \\

To prove the claim, we proceed by  induction on $j$\hskip0.5mm : 
This is obviously true for $j=1$ and $j=2$. 
Then by the induction hypothesis 
\begin{align*}Q_{r,j-1} &= r Q_{r,j-2} + Q_{r,j-3} < r (r+1)^{j-1} + (r+1)^{j-2} \\&= (r+1)^{j-2} (r (r+1) + 1) < (r+1)^{j-2} (r+1)^2 \\ &= (r+1)^{j}.\end{align*} 
In order to obtain an upper bound for the number $P(A,{\bf p}),$ 
it suffices to consider words over the largest allowed $s$-digit alphabet 
$A' : \{1 <\dots < s\}$,~ $b_1=t=l=1\ (< s)$, with Parikh vector 
${\bf p}' = (\underbrace{m,m,\ldots ,m}_{\mbox {$s$-times}}).$ 
Clearly
\vskip-0.1cm                                                                                             
\[
P(A,{\bf p}) \leq W_{\max}(A,{\bf p}) \leq W_{\max}(A',{\bf p'})
\]                                                                
\vskip-0.0cm 
and by (\ref{eqn2})  
\vskip-0.3cm
\begin{equation}\label{eqn3}
w_{\max}(A', {\bf p}') = s\cdot (s-1)^{m-1}\cdot (s-2) \cdots 1\cdot 1^{m-1}\cdots (s-2)^{m-1}\cdot (s-1)\cdot s^{m-1}.\end{equation}
\vskip0.2cm
By iteration of (\ref{ineq1}) 
applied to the decomposition in (\ref{eqn3})\ 
we obtain the inequalities
\vskip-0.4cm
\begin{align*}
W_{\max}(A', {\bf p}')&~ = K(w_{\max}(A', {\bf p}'))  \\ 
&~ < 2^{2s}~s\cdot (s-1)\cdots 3 \cdot 2~ \prod_{j=1}^s K(j^{m-1})\\ 
&~ = 2^{2s}~s!~ \prod_{j=1}^s Q_{j,m-1} \\
&~ < 2^{2s}~s!~ \prod_{j=1}^s (j+1)^m   \\ 
&~ =2^{2s}~s! ((s+1)!)^m 
\end{align*}
\vskip-0.2cm 
and hence
\begin{equation}\label{ineq4}
P(A,{\bf p}) < 2^{2s}\ s!\ ((s+1)!)^m.
\end{equation}
\vskip0.1cm
For each $s\geq 2$ we define $m_0 = m_0(s)$ to be the smallest positive integer such that
\[2^{2s}\ s! \leq \left (\frac{100}{99}\right)^{m_0}.\] 
\vskip-0.4cm 
Then 
\begin{equation}\label{ineq5}
P(A,{\bf p})  < \left(\frac{100}{99}~((s+1)!)\right)^m\ \text{\ for all\ ~} m\geq m_0(s).                                                                                                                           
\end{equation}                                                   

\vskip0.2cm  
On the other hand, we have the following lower bound for the number of different words in $\mathcal{X}(w)$: 
\vskip-0.15cm  
\begin{equation}\label{ineq6}
N(A,{\bf p})\geq\frac{((s-l+t)\hskip0.9mm m)!}{2(m!)^{s-l+t}}
\end{equation}                                                   

Based on the condition (\ref{ineq7})\ 
below, we will later  make a choice of $s = s'(t,l)$ depending on the parameters $t,l$. We apply the estimates provided by Sterling's formula to the factorial terms occurring in relations (\ref{ineq5}) and (\ref{ineq6})\ 
to obtain 
\begin{equation*}
(P(A,{\bf p}))^{1/m} < \frac{100}{99} (s+1)! < 
\frac{100}{99}\frac{12}{11}~{\textrm e}^{-(s+1)}(s+1)^{s+1}\sqrt{2\pi (s+1)}.
\end{equation*}
\vskip-0.35cm
\begin{equation*}
(((s-l+t)\hskip0.45mm m)!)^{1/m}
 > {\textrm e}^{-(s-l+t)} ((s-l+t)\hskip0.5mm m)^{s-l+t} 
{\sqrt{2 \pi (s-l+t)\hskip0.5mm m}}^{\hskip0.5mm 1/m}. 
\end{equation*} 
\vskip-0.6cm
\begin{equation*}
\left(2\hskip0.9mm (m!)^{s-l+t}\right)^{1/m} < {\textrm e}^{-(s-l+t)} m^{s-l+t} 
\left(2\ \frac{12}{11}\ \left(\sqrt{2\pi m}\right)^{s-l+t}\right)^{\hskip-0.4mm 1/m}.
\end{equation*} 
\vskip0.0cm
When we put the right hand sides of the last two inequalities together, the terms ${\textrm e}^{s-l+t}$ and $m^{s-l+t}$ cancel out, 
and if we keep the parameters $t, l$ fixed for the moment, 
the terms of the form $\sqrt{\ \cdot\ }^{~1/m}$ tend to 1 as $m\rightarrow\infty$. 
Letting $m\rightarrow\infty$ we get
\begin{align*}
\lim_{m\to\infty}\left(\frac{N(A,{\bf p})}{P(A,{\bf p})}\right)^{1/m}&\geq \frac{99}{100} \frac{11}{12} \frac{{\textrm e}^{s+1}(s-l+t)^{s-l+t}}{\sqrt{2\pi (s+1)}(s+1)^{s+1}}=\\
&=\frac{363}{400} \frac{{\textrm e}^{s+1}(s+1-l+t-1)^{~s+1-l+t-1}}
{\sqrt{2\pi(s+1)}(s+1)^{s+1}}\\
&=\frac{363}{400} \frac{{\textrm e}^{s+1}}{\sqrt{2\pi(s+1)}~(s-l+t)^{l-t+1}} 
\left(1-\frac{l-t+1}{s+1}\right)^{s+1}.
\end{align*} 
\vskip0.20cm
For fixed $t,l$\hskip2.5mm ($l-t\geq 1$)\hskip2.5mm the function $f(t,l,s) =
 \left(1-\frac{l-t+1}{s+1}\right)^{s+1}$ in the variable $s$ is strictly increasing on the interval $[l-t+1,\infty)$ with $f(t,l,s)\nearrow {\textrm e}^{-(l-t)-1}$ 
as $s\rightarrow \infty$.
We define $s_0$ to be the lowest integer such that $f(t,l,s_0)\geq 
\frac{1}{2}\ {\textrm e}^{-(l-t)-1}$. 
Then 

\parindent-0.0cm
\vskip-0.50cm
\begin{align*}
\lim_{m\to\infty}\left(\frac{N(A,{\bf p})}{P(A,{\bf p})}\right)^{1/m}\hskip-0.1cm
\geq\frac{363}{400} \frac{{\textrm e}^{s+1}}{\sqrt{2\pi(s+1)}~(s-l+t)^{l-t+1}}\ 
\frac{1}{2}\ {\textrm e}^{-(l-t)-1}\ =: H(t,l,s)
\end{align*}
\vskip-0.0cm
for all\ $s\geq s_0$. Obviously there exists some sufficiently large $s' = s'(t,l)\geq s_0$ 
such that 
\vskip-0.2cm
\begin{equation}\label{ineq7}
H(t,l,s') > 1.
\end{equation}                                                             
\vskip-0.3cm
Therefore the right hand side of
\vskip-0.4cm
\begin{equation}\label{ineq8}
\left( \frac{N(A,{\bf p})}{P(A,{\bf p})}\right) > (H(t,l,s'))^m                                 
\end{equation}                                                             
can be made arbitrarily large by letting $m\rightarrow \infty$. 
We call an $(s'-l+t)$-digit alphabet $A = \{(1\leq)~b_1<\dots < b_t <\dots < s'\}$
~{\it admissible} if $s' = s'(t,l)$ fulfills condition (\ref{ineq7})\ 
We consider the word $u(A,{\bf p}_1) = (b_1)^{m_1}\cdots (b_t)^{m_1}(l+1)^{m_1} \cdots (s')^{m_1}$ 
of length 
$n = (s'-l+t)\hskip0.5mm m_1$ 
with Parikh vector 
${\bf p}_1=((m_1)^{s'-l+t})$ 
where we choose $m_1\geq m_0$ such that $\left( \frac{N(A,{\bf p}_1)}{P(A,{\bf p}_1)}\right) > (H(t,l,s'))^{m_1}$.
The multi-set $\mathcal{X}_1=\mathcal{X}(A,{\bf p}_1)$ is made up of the $N(A,{\bf p}_1)=\#\mathcal{X}_1$ 
permuted arrangements of $u$. 
There exists at least one word $w_1\in\mathcal{X}_1$ with multiplicity $\mu\geq 2$ because otherwise we would have
$N(A,{\bf p}_1)=P(A,{\bf p}_1)$ which contradicts (\ref{ineq8})\  
with $m = m_1$. 
Let $\mu_1\ (\geq 2)$ be the maximal multiplicity attained by words $w\in\mathcal{X}_1$. 
Next choose $m_2 > m_1(s')$ such that $H(t,l,s')^{m_2} > \mu_1$. 
We claim that at least one word $w_2$ from $\mathcal{X}_2=\mathcal{X}(A,{\bf p}_2),
\ {\bf p}_2=((m_2)^{s'-l+t})$ has multiplicity $\mu > \mu_1$. 
Otherwise we would have $N(A,{\bf p}_2) \leq \mu_1 P(A,{\bf p}_2)$ which contradicts (\ref{ineq8})\  
with $m = m_2$. 
Next let $\mu_2\ (\geq \mu_1)$ be the maximal multiplicity attained by words $w\in\mathcal{X}_2$. 
Proceeding with this construction step by step we end up with a sequence of words $w_k$ on $A$ 
with multiplicities $\mu_k\to\infty$ as $k\to\infty$. The construction can be carried out for infinitely 
many different admissible alphabets.
This completes the proof of the Theorem.
\end{proof}
\medskip
The question remains largely unsolved in the case of semi-regular continuants though it seems certain that the behavior is quite similar to the regular case. 

\bigskip
There is some evidence supporting the following                     

\medskip
{\bf Conjecture.} 
Given any ordered alphabet $A =\{a_1<\dots <a_s\}~(a_j\in\mathbb{N},\ s\geq 2)$, let $\mu\geq 2$ be a positive integer. Then there exist infinitely many words on $A$ whose multiplicity is precisely $\mu$. 
The problem appears to require a difficult investigation into the values of continuants. 
Most likely our theorem and the conjecture also hold for continuants of semi-regular continued fractions. Unfortunately no higher-dimensional analogue of the theorem is available at present for $s\geq 4$ due to the fact that very little is known about the maximizing arrangements $w'_{\max}$ for $s\geq 4$ (see \cite{Zam}).


\bigskip
\hrulefill

{\small Gerhard Ramharter\\
Dept. of Convex and Discrete Geometry\\
Technische Universit\"at Wien\\
Wiedner Hauptstrasse 8-10\\
A-1040 Vienna, Austria\\
e-mail: Gerhard Ramharter $<$gerhard.ramharter@tuwien.ac.at$>$

\bigskip
Luca Q. Zamboni\\
Institut Camille Jordan\\
Universit\'e Claude Bernard Lyon 1\\
43 boulevard du 11 novembre 1918\\
F-69622 Villeurbanne Cedex\newline
e-mail: zamboni $<$zamboni@math.univ-lyon1.fr$>$}
\end{document}